\input amstex
\input amsppt.sty
\magnification=\magstep1
\hsize=30truecc
\vsize=22.2truecm
\baselineskip=16truept
\TagsOnRight \pageno=1 \nologo

\def\N{\Bbb N}

\def\l{\left}
\def\r{\right}
\def\bg{\bigg}
\def\({\bg(}
\def\[{\bg\lfloor}
\def\){\bg)}
\def\]{\bg\rfloor}
\def\t{\text}
\def\f{\frac}

\def\bi{\binom}
\def\eq{\equiv}

\def\ls{\leqslant}
\def\gs{\geqslant}
\def\mo{\roman{mod}}

\def\al{\alpha}
\def\da{\delta}

\def\jacob #1#2{\left(\frac{#1}{#2}\right)}

\def\pmod #1{\ (\roman{mod}\ #1)}

\def\Proof{\noindent{\it Proof}}

\def\Remark{\medskip\noindent{\it  Remark}}

\def\Ack{\medskip\noindent {\bf Acknowledgments}}
\hbox {Sci. China Math. 57(2014), no.\,10, 2091--2102.}
\bigskip
\topmatter
\title Proof of three conjectures on congruences\endtitle
\author Hao Pan and Zhi-Wei Sun\endauthor
\leftheadtext{Hao Pan and Zhi-Wei Sun}
\affil Department of Mathematics, Nanjing University\\
 Nanjing 210093, People's Republic of China
\\{\tt haopan1979\@gmail.com, \ zwsun\@nju.edu.cn}
\endaffil
\abstract In this paper we prove three conjectures on congruences
involving central binomial coefficients or Lucas sequences. Let
$p$ be an odd prime and let $a$ be a positive integer. We show that
if $p\equiv 1\pmod{4}$ or $a>1$ then
$$
\sum_{k=0}^{\lfloor\f34p^a\rfloor}\bi{-1/2}k\equiv\jacob{2}{p^a}\ \pmod{p^2},$$
where $(-)$ denotes the Jacobi symbol.
This confirms a conjecture of the second author. We also confirm a
conjecture of R. Tauraso by showing that
$$\sum_{k=1}^{p-1}\f{L_k}{k^2}\eq0\ \pmod{p}\quad \t{provided}\ p>5,$$ where the Lucas numbers
$L_0,L_1,L_2,\ldots$ are defined by $L_0=2,\ L_1=1$ and
$L_{n+1}=L_n+L_{n-1}\ (n=1,2,3,\ldots)$. Our third theorem
states that if $p\not=5$ then we can determine $F_{p^a-(\f{p^a}5)}$ mod $p^3$  in the following way:
$$\sum_{k=0}^{p^a-1}(-1)^k\bi{2k}k\eq\l(\f{p^a}5\r)\l(1-2F_{p^a-(\f{p^a}5)}\r)\ \pmod{p^3},$$
which appeared as a conjecture in a paper of Sun and Tauraso in 2010.
\endabstract
\thanks 2010 {\it Mathematics Subject Classification}.\,Primary 11B65, 11A07;
Secondary 05A10, 11B39.
\newline\indent The second author is the corresponding author.
\endthanks
\endtopmatter
\document

\heading{1. Introduction}\endheading

In this paper we aim to prove three conjectures on congruences.

Our first theorem confirms a conjecture raised by the second author
[S11, Conjecture 1.2].

\proclaim{Theorem 1.1} Let $p$ be an odd prime and let $a$ be a
positive integer. If $p\equiv 1\pmod{4}$ or $a>1$, then we have
$$\sum_{k=0}^{\lfloor\f34p^a\rfloor}\bi{-1/2}k\equiv\jacob{2}{p^a}\pmod{p^2},
\tag1.1$$
where $(-)$ denotes the Jacobi symbol.
\endproclaim

Our second theorem confirms a nice conjecture of R. Tauraso [T], and it presents a congruence involving Lucas numbers
which is similar to the well-known Wolstenholme congruence $\sum_{k=1}^{p-1}1/k^2\eq0\ (\mo\ p)$ with $p>3$ prime (cf. [Wo]).

\proclaim{Theorem 1.2} Let $p>5$ be a prime. Then
$$\sum_{k=1}^{p-1}\f{L_k}{k^2}\eq0\ \pmod{p},\tag1.2$$
where the Lucas numbers $L_0,L_1,L_2,\ldots$ are defined by $L_0=2,\
L_1=1$ and $L_{n+1}=L_n+L_{n-1}\ (n=1,2,3,\ldots)$.
\endproclaim

The Fibonacci sequence $\{F_n\}_{n\gs0}$ is given by $F_0=0,\ F_1=1$
and $F_{n+1}=F_n+F_{n-1}\ (n=1,2,3,\ldots)$. It is well-known that
$p\mid F_{p-(\f p5)}$ for any odd prime $p$, and the Fibonacci quotient $F_{p-(\f p5)}/p$ modulo $p$ is closely related to
fundamental units of real quadratic fields (cf. H. C. Williams [W])
and Vandiver's conjecture about class numbers of real cyclotomic fields
(cf. S. Jakubec [J]). Our following theorem
determines $F_{p^a-(\f {p^a}5)}$ mod $p^3$ for any $a=1,2,3,\ldots$, and the result appeared as
Conjecture 1.1 of Sun and Tauraso [ST].

\proclaim{Theorem 1.3} Let $p\not= 2,5$ be a prime and let $a$ be a positive integer. Then
we have
$$\sum_{k=0}^{p^a-1}(-1)^k\bi{2k}k\eq\l(\f{p^a}5\r)(1-2F_{p^a-(\f{p^a}5)})\ \ (\mo\ p^3).\tag1.3$$
\endproclaim

Note that (1.3) modulo $p$ is [ST, (1.7)] with $d=0$, and (1.3) modulo $p^2$ was given in [S10, Corollary 1.1].

Those primes $p>5$ satisfying the congruence $F_{p-(\f p5)}\eq0\
(\mo\ p^2)$ are called Wall-Sun-Sun primes (cf. [CP, p.\,32]
and [SS]). It is known that there are no Wall-Sun-Sun primes below
$4.5\times 10^{16}$ (cf. [P]).

To understand our proofs of Theorems 1.1-1.3 one needs some
basic knowledge of Lucas sequences.

Let $A$ and $B$ be two integers. The Lucas sequence $u_n=u_n(A,B)$
$(n\in\N=\{0,1,2,\ldots\})$ and its companion $v_n=v_n(A,B)\
(n\in\N)$ are defined by
$$u_0=0,\ u_1=1,\ u_{n+1}=Au_{n}-Bu_{n-1}\ \ \t{for}\ n=1,2,3,\ldots,$$
and
$$v_0=2,\ v_1=A,\ v_{n+1}=Av_{n}-Bv_{n-1}\ \ \t{for}\ n=1,2,3,\ldots.$$
Let $\Delta=A^2-4B$. The characteristic equation $x^2-Ax+B=0$ has two roots
$$\al=\f{A+\sqrt{\Delta}}2\ \ \ \t{and}\ \ \
\beta=\f{A-\sqrt{\Delta}}2,$$ which are both algebraic integers.
It is well known that
$$(\al-\beta)u_n=\al^n-\beta^n\ \ \t{and}\ \ v_n=\al^n+\beta^n\ \
\t{for all}\ n=0,1,2,\ldots.$$ For an odd prime $p$ and a positive integer $a$, clearly
$$v_{p^a}\eq(\al+\beta)^{p^a}=A^{p^a}\eq A\ (\mo\ p)$$
and
$$\Delta u_{p^a}=(\al-\beta)(\al^{p^a}-\beta^{p^a})\eq(\al-\beta)^{p^a+1}
=\Delta^{(p^a+1)/2}\eq\Delta\l(\f{\Delta}{p^a}\r)\ (\mo\ p).$$ It is
also known that $p^a\mid u_{p^a-(\f{\Delta}{p^a})}$ provided that
$p\nmid B$ (see, e.g., [S10, Lemma 2.3]).

Note that $F_n=u_n(1,-1)\ (n\in\N)$ and $L_n=v_n(1,-1)\ (n\in\N)$
are familiar Fibonacci numbers and Lucas numbers respectively. The
Pell sequence and its companion are given by $P_n=u_n(2,-1)\
(n\in\N)$ and $Q_n=v_n(2,-1)\ (n\in\N)$ respectively.

We will show Theorems 1.1-1.3 in Sections 2-4 respectively. In the proofs of Theorems 1.2 and 1.3,
we employ the useful technique of A. Granville [Gr] and deal with congruences in the ring of algebraic integers.

\heading{2. Proof of Theorem 1.1}\endheading

\proclaim{Lemma 2.1} Let $p$ be an odd prime and let $a$ be a positive integer. Then
$$P_{p^a-\jacob{2}{p^a}}Q_{p^a-\jacob{2}{p^a}}\equiv\jacob{2}{p^a}\f{Q_{p^a}-2}2 \pmod{p^2}.\tag2.1$$
\endproclaim
\Proof. Recall that $P_n=u_n(2,-1)$ and $Q_n=v_n(2,-1)$ for all $n\in\N$. So
$$P_{p^a}\equiv\jacob2{p^a}\pmod{p}\ \ \ \t{and}\ \ \
Q_{p^a}\equiv2\pmod{p}.$$
Since $Q_{n-1}=4P_n-Q_n$ and $Q_{n+1}=4P_n+Q_n$ for $n=1,2,3,\ldots$, we have
$$\jacob{2}{p^a}Q_{p^a-\jacob{2}{p^a}}=4\jacob{2}{p^a}P_{p^a}-Q_{p^a}\equiv 2\pmod{p}.$$
Similarly,
$$P_{p^a-(\f2{p^a})}=\f{Q_{p^a}}2-\jacob{2}{p^a}P_{p^a}\equiv0\pmod{p}.$$
Clearly $\al=1+\sqrt2$ and $\beta=1-\sqrt2$ are the two roots of the equation $x^2-2x-1=0$.
Thus
$$Q_n^2-8P_n^2=(\al^n+\beta^n)^2-(\al^n-\beta^n)^2=4(\al\beta)^n=4(-1)^n$$
for all $n\in\N$.
Therefore
$$Q_{p^a-\jacob{2}{p^a}}^2-4=8P_{p^a-\jacob2{p^a}}^2\equiv0\pmod{p^2}$$
and hence
$$4\jacob{2}{p^a}P_{p^a}-Q_{p^a}=\jacob{2}{p^a}Q_{p^a-\jacob{2}{p^a}}\equiv2\pmod{p^2}.
$$
It follows that
$$
\align
&\jacob{2}{p^a}P_{p^a-\jacob{2}{p^a}}Q_{p^a-\jacob{2}{p^a}}
\\\eq&2P_{p^a-\jacob{2}{p^a}}
= Q_{p^a}-2\jacob{2}{p^a}P_{p^a}\equiv\frac{Q_{p^a}}2-1\pmod{p^2}.
\endalign$$
This proves (2.1). \qed

\proclaim{Lemma 2.2} Let $p$ be an odd prime and let $a$ be a positive integer.
Suppose that $p\eq1\ (\mo\ 4)$ or $a>1$. Then
$$p^a\sum_{0\ls k<\lfloor p^a/4\rfloor}\f1{\bi{(p^a-3)/2}k}\eq\l(\f2{p^a}\r)\f{Q_{p^a}-2}4\ (\mo\ p^2).\tag2.2$$
\endproclaim
\Proof. If $p^a\eq1\ (\mo\ 4)$, then $(p^a-3)/2$ is odd and hence
$$\sum_{k=0}^{(p^a-1)/4-1}\frac{1}{\binom{(p^a-3)/2}{k}}=\frac12\sum_{k=0}^{(p^a-3)/2}\frac{1}{\binom{(p^a-3)/2}{k}}.
$$
If $p^a\eq3\ (\mo\ 4)$, then $a\in\{3,5,\ldots\}$ and
$$
\sum_{k=0}^{(p^a-3)/4-1}\frac{1}{\binom{(p^a-3)/2}{k}}=\frac12\sum_{k=0}^{(p^a-3)/2}\frac{1}{\binom{(p^a-3)/2}{k}}-\frac12\cdot\frac{1}{\binom{(p^a-3)/2}{(p^a-3)/4}}.
$$
In the case $p^a\eq3\ (\mo\ 4)$,  as the fractional parts of $(p^a-3)/(2p)$ and $(p^a-3)/(4p)$ are $(p-3)/(2p)$ and $(p-3)/(4p)$ respectively,
we have
$$\l\lfloor\f{(p^a-3)/2}p\r\rfloor=2\l\lfloor\f{(p^a-3)/4}p\r\rfloor$$
and hence
$$\nu_p\bigg(\binom{(p^a-3)/2}{(p^a-3)/4}\bigg)=\sum_{j=1}^{a-1}\bigg(\left\lfloor{\frac{(p^a-3)/2}{p^j}}\right\rfloor-2\left\lfloor{\frac{(p^a-3)/4}{p^j}}\right\rfloor\bigg)<a-1,
$$
where $\nu_p(x)$ denotes the $p$-adic valuation of an integer $x$. (It is well known that
$\nu_p(n!)=\sum_{j=1}^\infty\lfloor n/p^j\rfloor$ for any $n\in\N$.) Therefore,
$$p^a\sum_{0\ls k<\lfloor p^a/4\rfloor}\f1{\bi{(p^a-3)/2}k}\eq\f{p^a}2\sum_{k=0}^{(p^a-3)/2}\f1{\bi{(p^a-3)/2}k}\ (\mo\ p^2).$$

Applying the known identity
$$
\sum_{k=0}^n\frac{x^k}{\binom{n}{k}}=(n+1)\l(\frac{x}{1+x}\r)^{n+1}\,\sum_{k=1}^{n+1}\frac{1+x^k}{k(1+x)}\bigg(\frac{1+x}x\bigg)^k
$$
(cf. [G, (2.4)]) with $x=1$ and $n=(p^a-3)/2$, we get
$$\align
\frac{p^a}2\sum_{k=0}^{(p^a-3)/2}\frac{1}{\binom{(p^a-3)/2}{k}}
=&\frac{p^a(p^a-1)}{2^{(p^a+3)/{2}}}\sum_{k=1}^{(p^a-1)/2}\frac{2^k}k
\\\equiv&-\jacob{2}{p^a}\frac{p^a}{4}\sum_{k=1}^{(p^a-1)/2}\frac{2^k}k
\pmod{p^2}.
\endalign$$
Since
$$\bi{p^a}j=\f{p^a}j\prod_{0<i<j}\f{p^a-i}i\eq\f{p^a}j(-1)^{j-1}\ (\mo\ p^2)$$
for all $j=1,\ldots,p^a-1$, we have
$$\align
p^{a}\sum_{k=1}^{(p^a-1)/2}\frac{2^k}k\equiv&
-2\sum_{k=1}^{(p^a-1)/2}\binom{p^a}{2k}2^k
=-\sum_{k=1}^{p^a}\binom{p^a}{k}(\sqrt{2}^k+(-\sqrt{2})^k)\\
=&-(1+\sqrt2)^{p^a}-(1-\sqrt2)^{p^a}+2=-Q_{p^a}+2\pmod{p^2}.
\endalign$$

Combining the above we immediately get (2.2). \qed

\medskip
\noindent{\it Proof of Theorem 1.1}. Clearly,
$$\bi{-1/2}k=\f{\bi{2k}k}{(-4)^k}\qquad \ \t{for all}\ k=0,1,2,\ldots.$$
Choose $\delta\in\{1,3\}$ such that $p^a\equiv\delta\pmod{4}$. Then
$$\sum_{k=0}^{\lfloor\f 34p^a\rfloor}\bi{-1/2}k=
\sum_{k=0}^{p^a-1}\frac{\bi{2k}k}{(-4)^k}-\sum_{k=(3p^a+\delta)/4}^{p^a-1}\frac{\bi{2k}k}{(-4)^k}.
$$
By Sun [S10, Theorem 1.1 and Lemma 2.3],
$$
\sum_{k=0}^{p^a-1}\frac{\bi{2k}k}{(-4)^k}\equiv \jacob{2}{p^a}+u_{p^a-\jacob{2}{p^a}}(-6,1)\pmod{p^2}.
$$
Hence we only need to prove the following congruence:
$$\sum_{k=(3p^a+\delta)/4}^{p^a-1}\frac{\bi{2k}k}{(-4)^k}\equiv u_{p^a-\jacob{2}{p^a}}(-6,1)\pmod{p^2}\tag2.3$$
provided that $p\eq1\ (\mo\ 4)$ or $a>1$.

Let $k$ and $l$ be positive integers with $k+l=p^a$ and $0<l<p^a/2$. Then
$$\f{\bi{2k}k}{\bi{2p^a-2}{p^a-1}}=\f{(2p^a-2l)!}{(2p^a-2)!}\(\f{(p^a-1)!}{(p^a-l)!}\)^2=\f{\prod_{0<i<l}(p^a-i)^2}{\prod_{1<j<2l}(2p^a-j)}$$
and hence
$$\f{\bi{2k}k}{\bi{2p^a-2}{p^a-1}}\cdot\f{(2l-1)!}{(l-1)!^2}=\f{\prod_{0<i<l}(1-p^a/i)^2}{\prod_{1<j<2l}(1-2p^a/j)}\eq1\
(\mo\ p).$$ Note that
$$
\binom{2p^a-2}{p^a-1}=p^a\prod_{j=2}^{p^a-1}\frac{2p^a-j}{j}\equiv-p^a\pmod{p^{a+1}}
$$
and
$$\bi{2k}k=\bi{p^a+(2k-p^a)}{0p^a+k}\eq\bi{2k-p^a}k=0\ (\mo\ p)$$
by Lucas' theorem. So we have
$$\f l2\bi{2l}l=\f{(2l-1)!}{(l-1)!^2}\not\eq0\ (\mo\ p^a)$$
and
$$\bi{2k}k\eq-p^a\f{(l-1)!^2}{(2l-1)!}=-\f{2p^a}{l\bi{2l}l}\ (\mo\ p^2).$$

In view of the above,
$$\sum_{k=(3p^a+\delta)/4}^{p^a-1}\frac{\bi{2k}k}{(-4)^k}
\equiv\frac{-2p^a}{(-4)^{p^a}}\sum_{l=1}^{(p^a-\delta)/4}\frac{(-4)^l}{l\binom{2l}{l}}
\eq \f{p^a}2 \sum_{k=1}^{(p^a-\delta)/4}\frac{(-4)^k}{k\binom{2k}{k}} \pmod{p^2}.$$
For $k=1,\ldots,(p^a-1)/2$, clearly
$$\align
\frac{\binom{(p^a-1)/2}{k}}{\binom{2k}k/(-4)^k}=&\frac{\binom{(p^a-1)/2}{k}}{\binom{-1/2}{k}}=\prod_{j=0}^{k-1}\frac{(p^a-1)/2-j}{-1/2-j}\\
=&\prod_{j=0}^{k-1}\bigg(1-\frac{p^a}{2j+1}\bigg)\equiv1\pmod{p}
\endalign$$
and hence
$$\f{\bi{(p^a-3)/2}{k-1}}{k\bi{2k}k/(-4)^k}\eq\f2{p^a-1}\eq-2\ (\mo\ p).$$
Therefore
$$\frac{p^a}{2}\sum_{k=1}^{(p^a-\delta)/4}\frac{(-4)^k}{k\binom{2k}{k}}
\equiv-p^a\sum_{k=0}^{(p^a-\delta)/4-1}\frac{1}{\binom{(p^a-3)/2}{k}}\pmod{p^2}.$$
So far we have reduced (2.3) to the following congruence:
$$p^a\sum_{k=0}^{(p^a-\da)/4-1}\f1{\bi{(p^a-3)/2}k}\eq-u_{p^a-(\f2{p^a})}(-6,1)\ (\mo\ p^2).\tag2.4$$

 In view of (2.4) and Lemma 2.2, it suffices to show that
 $$u_{p^a-(\f2{p^a})}(-6,1)\eq-\l(\f2{p^a}\r)\f{Q_{p^a}-2}4\ (\mo\ p^2).$$
As $-3+2\sqrt2$ and $-3-2\sqrt2$ are the two roots of the equation $x^2+6x+1=0$, for any $n\in\N$ we have
$$\align u_n(-6,1)=&\f{(-3+2\sqrt2)^n-(-3-2\sqrt2)^n}{4\sqrt2}
\\=&\frac{(-1)^{n-1}}2\cdot\frac{(1+\sqrt{2})^{2n}-(1-\sqrt{2})^{2n}}{2\sqrt{2}}=\frac{(-1)^{n-1}}{2}P_{n}Q_n
\endalign$$
Therefore, with the help of Lemma 2.1 we finally obtain
$$u_{p^a-(\f2{p^a})}(-6,1)=-\f12P_{p^a-(\f2{p^a})}Q_{p^a-(\f2{p^a})}\eq-\l(\f2{p^a}\r)\f{Q_{p^a}-2}4\ (\mo\ p^2)$$
as desired.

The proof of Theorem 1.1 is now complete. \qed

\heading{3. Proof of Theorem 1.2}\endheading

\proclaim{Lemma 3.1} Let $p>3$ be a prime. Then we have the following congruence:
$$\l(\f{x^p+(1-x)^p-1}p\r)^2\eq-2\sum_{k=1}^{p-1}\f{(1-x)^k}{k^2}-2x^{2p}\sum_{k=1}^{p-1}\f{(1-x^{-1})^k}{k^2}\ (\mo\ p).
\tag3.1$$
\endproclaim
\Proof. (3.1) follows immediately if we combine (4) and (5) of Granville [Gr]. \qed

\proclaim{Proposition 3.2} Let $A$ and $B$ be nonzero integers, and let $\al$ and $\beta$ be the two roots of the equation $x^2-Ax+B=0$.
Let $p$ be an odd prime not dividing $AB$. Then
$$\l(\f{v_p(A,B)-A^p}p\r)^2\eq-2A^2\sum_{k=1}^{p-1}\f{\al^k}{A^kk^2}
-2\beta^{2p}\sum_{k=1}^{p-1}\f{\al^{2k}}{(-B)^kk^2}\ (\mo\ p),\tag3.2$$
and
$$\l(\f{v_p(A,B)-A^p}p\r)^2\eq-2A\al^p\sum_{k=1}^{p-1}\f{\al^k}{A^kk^2}
-2\beta^{2p}\sum_{k=1}^{p-1}\f{A^k\al^k}{B^kk^2}\ (\mo\ p).\tag3.3$$
\endproclaim
\Proof. By (3.1) and Fermat's little theorem,
$$\align&\f1{A^2}\l(\f{x^p+(A-x)^p-A^p}{p}\r)^2
\\\eq&\l(\f{(x/A)^p+(1-x/A)^p-1}p\r)^2
\\\eq&-2\sum_{k=1}^{p-1}\f{(1-x/A)^k}{k^2}-2\l(\f xA\r)^{2p}\sum_{k=1}^{p-1}\f{(1-A/x)^k}{k^2}\ (\mo\ p).
\endalign$$
Note that $v_p(A,B)=\beta^p+\al^p=\beta^p+(A-\beta)^p$ and $\al\beta=B$. So we have
$$\align\l(\f{v_p(A,B)-A^p}p\r)^2\eq&-2A^2\sum_{k=1}^{p-1}\f{(A-\beta)^k}{A^kk^2}
\\&-2\beta^{2p}\sum_{k=1}^{p-1}\f{(1-A\al/B)^k}{k^2}
\ (\mo\ p)\endalign$$
and hence (3.2) holds since $A\al-B=\al^2$.

On the other hand,
$$\al^p(A^p-v_p(A,B))=\al^p(A^p-\al^p-\beta^p)=(B+\al^2)^p+(-\al^2)^p-B^p$$
and hence
$$\align&\al^{2p}\l(\f{A^p-v_p(A,B)}p\r)^2
\\=&\l(\f{(-\al^2)^p+(B-(-\al^2))^p-B^p}p\r)^2
\\\eq&-2B^2\sum_{k=1}^{p-1}\f{(1-(-\al^2)/B)^k}{k^2}-2(-\al^2)^{2p}\sum_{k=1}^{p-1}\f{(1-B/(-\al^2))^k}{k^2}
\\=&-2B^2\sum_{k=1}^{p-1}\f{(A\al)^k}{B^kk^2}-2\al^{4p}\sum_{k=1}^{p-1}\f{(A\al)^k}{\al^{2k}k^2}
\\\eq&-2(\al\beta)^{2p}\sum_{k=1}^{p-1}\f{(A\al)^k}{B^kk^2}-2A\al^{3p}\sum_{k=1}^{p-1}\f{\al^{p-k}}{A^{p-k}(p-k)^2}\ (\mo\ p).
\endalign$$
Therefore (3.3) follows. \qed

\medskip
\noindent{\it Proof of Theorem 1.2}. Let $\al$ and $\beta$ be the two roots of the equation $x^2-x-1=0$.
Applying Proposition 3.2 with $A=1$ and $B=-1$, we get
$$\l(\f{L_p-1}p\r)^2\eq-2\sum_{k=1}^{p-1}\f{\al^k}{k^2}-2\beta^{2p}\sum_{k=1}^{p-1}\f{\al^{2k}}{k^2}\ (\mo\ p)\tag3.4$$
and
$$\l(\f{L_p-1}p\r)^2\eq-2\al^p\sum_{k=1}^{p-1}\f{\al^k}{k^2}-2\beta^{2p}\sum_{k=1}^{p-1}\f{(-\al)^k}{k^2}\ (\mo\ p).\tag3.5$$
Since
$$\align\sum_{k=1}^{p-1}\f{2\al^{2k}}{(2k)^2}=&\sum_{j=1}^{2p-1}(1+(-1)^j)\f{\al^{j}}{j^2}
=\sum_{k=1}^{p-1}\(\f{\al^k+(-\al)^k}{k^2}+\f{\al^{p+k}+(-\al)^{p+k}}{(p+k)^2}\)
\\\eq&(1+\al^p)\sum_{k=1}^{p-1}\f{\al^k}{k^2}+(1-\al^p)\sum_{k=1}^{p-1}\f{(-\al)^k}{k^2}\ (\mo\ p),
\endalign$$
(3.4) can be rewritten as
$$\aligned\l(\f{L_p-1}p\r)^2\eq&-2(1+2(1+\al^p)\beta^{2p})\sum_{k=1}^{p-1}\f{\al^k}{k^2}
\\&-4(1-\al^p)\beta^{2p}\sum_{k=1}^{p-1}\f{(-\al)^k}{k^2}
\ (\mo\ p).\endaligned\tag3.6$$
Multiplying (3.5) by $2(1-\al^p)$ and then subtracting it from (3.6) we obtain
$$\align(2\al^p-1)\l(\f{L_p-1}p\r)^2\eq&\l(4\al^p(1-\al^p)-2-4(1+\al^p)\beta^{2p}\r)\sum_{k=1}^{p-1}\f{\al^k}{k^2}
\\=&(4L_p-4L_{2p}-2)\sum_{k=1}^{p-1}\f{\al^k}{k^2}\ (\mo\ p).
\endalign$$
Now that $L_p\eq1\ (\mo\ p)$ and
$$L_{2p}=\al^{2p}+\beta^{2p}\eq(\al^2+\beta^2)^p=\l((\al+\beta)^2-2\al\beta\r)^p=3^p\eq3\ (\mo\ p),$$
we have
$$(2\al^p-1)\l(\f{L_p-1}p\r)^2\eq(4-4\times3-2)\sum_{k=1}^{p-1}\f{\al^k}{k^2}=-10\sum_{k=1}^{p-1}\f{\al^k}{k^2}\ (\mo\ p).$$
Similarly,
$$(2\beta^p-1)\l(\f{L_p-1}p\r)^2\eq-10\sum_{k=1}^{p-1}\f{\beta^k}{k^2}\ (\mo\ p).\tag3.7$$
As $2\al^p-1+(2\beta^p-1)=2L_p-2\eq0\ (\mo\ p)$, we finally obtain
$$\sum_{k=1}^{p-1}\f{L_k}{k^2}=\sum_{k=1}^{p-1}\f{\al^k+\beta^k}{k^2}\eq0\ (\mo\ p).$$

So far we have completed the proof of Theorem 1.2. \qed

\Remark\ 3.3. Let $p>5$ be a prime. In view of (3.7), we also have
$$\sum_{k=1}^{p-1}\f{F_k}{k^2}\eq-\f{2F_p}{10}\l(\f{L_p-1}p\r)^2\eq-\f15\l(\f p5\r)\l(\f{L_p-1}p\r)^2\ (\mo\ p).\tag3.8$$

\heading{4. Proof of Theorem 1.3}\endheading

We need a lemma which extends a congruence due to Granville [Gr, (6)].

\proclaim{Lemma 4.1} Let $p$ be an odd prime and let $a$ be a positive integer.
Then
$$\aligned &p\da_{p,3}+p^{a-1}\sum_{k=1}^{p^a-1}\f{(1-x)^k}k
\\\eq&\f{1-x^{p^a}-(1-x)^{p^a}}p
-p\(\sum_{k=1}^{p-1}\f{x^{k}}{k^2}\)^{p^{a-1}}
\ (\mo\ p^2),
\endaligned\tag4.1$$
where the Kronecker symbol $\da_{p,3}$ is $1$ or $0$ according as $p=3$ or not.
\endproclaim
\Proof. As observed by Granville [Gr], for any integer $n>1$ we have
$$\align\sum_{k=1}^{n-1}\f{(1-x)^k-1}k=&\sum_{k=1}^{n-1}\f1k\sum_{j=1}^k\bi kj(-x)^j=\sum_{j=1}^{n-1}\f{(-x)^j}j\sum_{k=j}^{n-1}\bi{k-1}{j-1}
\\=&\sum_{j=1}^{n-1}\f{(-x)^j}j\bi{n-1}j\ \ (\t{by [G, (1.52)]})
\\=&\sum_{j=1}^{n-1}\f{(-x)^j}j\(\f nj\bi{n-1}{j-1}-\f jn\bi nj\)
\\=&n\sum_{j=1}^{n-1}\bi{n-1}{j-1}\f{(-x)^j}{j^2}-\f{(1-x)^n-(-x)^n-1}n.
\endalign$$
Note that
$$\align2p^{a-1}\sum_{k=1}^{p^a-1}\f1k=&p^{a-1}\sum_{k=1}^{p^a-1}\l(\f1k+\f1{p^a-k}\r)=\sum_{k=1}^{p^a-1}\f{p^a}k\cdot\f{p^{a-1}}{p^a-k}
\\\eq&\sum_{j=1}^{p-1}\f{p^a}{p^{a-1}j}\cdot\f{p^{a-1}}{p^a-p^{a-1}j}\eq-p\sum_{j=1}^{p-1}\f1{j^2}\eq p\da_{p,3}\ (\mo\ p^2).
\endalign$$
(As $\sum_{j=1}^{p-1}1/j^2\eq\sum_{k=1}^{p-1}1/(2k)^2\ (\mo\ p)$, we have $\sum_{j=1}^{p-1}1/j^2\eq0\ (\mo\ p)$ if $p>3$.)
Hence
$$p^{a-1}\sum_{k=1}^{p^a-1}\f1k\eq-p\da_{p,3}\ (\mo\ p^2).$$
Note also that
$$\bi{p^a-1}{j-1}=\prod_{0<i<j}\f{p^a-i}i\eq(-1)^{j-1}\ (\mo\ p).$$
Therefore
$$\align &p^{a-1}\sum_{k=1}^{p^a-1}\f{(1-x)^k}k+p\da_{p,3}
\\\eq &p^{a-1}\sum_{k=1}^{p^a-1}\f{(1-x)^k-1}k
\\\eq&-p^{2a-1}\sum_{j=1}^{p^a-1}\f{x^j}{j^2}-\f{x^{p^a}+(1-x)^{p^a}-1}{p}\ (\mo\ p^2)
\endalign$$
To complete the proof, it suffices to note that
$$p^{2a-1}\sum_{j=1}^{p^a-1}\f{x^j}{j^2}\eq p^{2a-1}\sum_{k=1}^{p-1}\f{x^{p^{a-1}k}}{(p^{a-1}k)^2}
=p\sum_{k=1}^{p-1}\f{x^{p^{a-1}k}}{k^2}
\ (\mo\ p^2)$$
and
$$\(\sum_{k=1}^{p-1}\f{x^k}{k^2}\)^{p^{a-1}}
\eq\sum_{k=1}^{p-1}\l(\f{x^k}{k^2}\r)^{p^{a-1}}\eq\sum_{k=1}^{p-1}\f{x^{p^{a-1}k}}{k^2}
\ (\mo\ p).$$
This concludes the proof. \qed

\proclaim{Proposition 4.2} Let $p\not=2,5$ be a prime and let $a$ be a positive integer. Then
$$p^{a-1}\sum_{k=1}^{p^a-1}\f{F_{2(p^a-k)}}{k}\eq\f{F_{2p^a}-F_{p^a}}p+\f p{10}\l(\f {p^a}5\r)\l(\f{L_p-1}p\r)^2\ (\mo\ p^2).\tag4.2$$
\endproclaim
\Proof. Let $\al$ and $\beta$ be the two roots of the equation $x^2-x-1=0$.
Clearly $\al+\beta=1$ and $\al\beta=-1$.
Set
$$g(x)=p^{a-1}\sum_{k=1}^{p^a-1}\f{x^k}k,\ q(x)=\f{x^{p^a}+(1-x)^{p^a}-1}p \ \t{and}\ G(x)=\sum_{k=1}^{p-1}\f{x^k}{k^2}.$$
By Lemma 4.1 we have
$$p\da_{p,3}+g(1-x)\eq-q(x)-pG(x)^{p^{a-1}}\ (\mo\ p^2).$$
In view of (3.7),
$$G(\beta)\eq\f{1-2\beta^p}{10}\l(\f{L_p-1}p\r)^2-\da_{p,3}\ (\mo\ p)$$
(this can be checked directly when $p=3$).
Hence
$$\align G(-\al)=&G(\beta^{-1})=\f1{\beta^p}\sum_{k=1}^{p-1}\f{\beta^{p-k}}{k^2}
\\\eq&-\al^pG(\beta)
\eq\da_{p,3}\al^p-\f{2+\al^p}{10}\l(\f{L_p-1}p\r)^2\ (\mo\ p).
\endalign$$
Note that
$$q(-\al)=\f{(-\al)^{p^a}+(1+\al)^{p^a}-1}p=\f{\al^{p^a}(\al^{p^a}+\beta^{p^a}-1)}p=\al^{p^a}\f{L_{p^a}-1}p.$$
Applying Lemma 4.1 we get
$$\align &g(\al^2)=g(1-(-\al))
\\\eq&-p\da_{p,3}-pG(-\al)^{p^{a-1}}-q(-\al)
\\\eq&-p\da_{p,3}-p\da_{p,3}\al^{p^a}+p\(\f{\al^p+2}{10}\)^{p^{a-1}}\l(\f{L_p-1}p\r)^{2p^{a-1}}
 -\al^{p^a}\f{L_{p^a}-1}{p}
\\\eq&-p\da_{p,3}(1+\al)^{p^a}+p\f{\al^{p^a}+2}{10}\l(\f{L_p-1}p\r)^2-\al^{p^a}\f{L_{p^a}-1}{p}\ (\mo\ p^2)
\endalign$$
and hence
$$\align&\beta^{2p^a}g(\al^2)+p\da_{p,3}(\al\beta)^{2p^a}
\\\eq&p\f{2\beta^{2p^a}-\beta^{p^a}}{10}\l(\f{L_p-1}p\r)^2+\beta^{p^a}\f{L_{p^a}-1}p\ (\mo\ p^2).
\endalign$$
As $\beta^{2p^a}=(1+\beta)^{p^a}\eq 1+\beta^{p^a}\ (\mo\ p)$, we have
$$\beta^{2p^a}g(\al^2)\eq -p\da_{p,3}+p\f{2+\beta^{p^a}}{10}\l(\f{L_p-1}p\r)^2+\beta^{p^a}\f{L_{p^a}-1}p\ (\mo\ p^2).$$
Similarly,
$$\al^{2p^a}g(\beta^2)\eq -p\da_{p,3}+p\f{2+\al^{p^a}}{10}\l(\f{L_p-1}p\r)^2+\al^{p^a}\f{L_{p^a}-1}p\ (\mo\ p^2).$$

Observe that
$$\align \sum_{k=1}^{p^a-1}\f{F_{2(p^a-k)}}k
=&\sum_{k=1}^{p^a-1}\f{\al^{2p^a-2k}-\beta^{2p^a-2k}}{(\al-\beta)k}
\\=&\f1{\al-\beta}\(\al^{2p^a}\sum_{k=1}^{p^a-1}\f{\beta^{2k}}k-\beta^{2p^a}\sum_{k=1}^{p^a-1}\f{\al^{2k}}k\).
\endalign$$
So, by the above, we have
$$\align p^{a-1}\sum_{k=1}^{p^a-1}\f{F_{2(p^a-k)}}{k}=&\f{\al^{2p^a}g(\beta^2)-\beta^{2p^a}g(\al^2)}{\al-\beta}
\\\eq&p\f{F_{p^a}}{10}\l(\f{L_p-1}p\r)^2+F_{p^a}\f{L_{p^a}-1}p
\\\eq&\f p{10}\l(\f{p^a}5\r)\l(\f{L_p-1}p\r)^2+\f{F_{2p^a}-F_{p^a}}p
\ (\mo\ p^2).
\endalign$$
This concludes the proof. \qed

\proclaim{Lemma 4.3} Let $p\not=2,5$ be a prime and let $a$ be a positive integer. Then
$$\l(\f{p^a}5\r)(2F_{p^a}-F_{2p^a})+\f{(L_p-1)^2}5\eq 1-2F_{p^a-(\f{p^a}5)}\ (\mo\ p^3).\tag4.3$$
\endproclaim
\Proof. Note that
$$(L_{p^a}-1)^2=p^2\l(\f{L_{p^a}-1}p\r)^2\eq p^2\l(\f{L_p-1}p\r)^2=(L_p-1)^2\ (\mo\ p^3)$$
since $L_{p^a}\eq L_p\ (\mo\ p^2)$ by [S10, (2.4)]. Also,
$$L_{p^a}=F_{p^a}+2F_{p^a-1}=2F_{p^a+1}-F_{p^a}=2F_{p^a-(\f {p^a}5)}+\l(\f{p^a}5\r)F_{p^a}.$$
Thus
$$\align &1-2F_{p^a-(\f{p^a}5)}-\l(\f{p^a}5\r)(2F_{p^a}-F_{2p^a})
\\=&1-L_{p^a}+\l(\f{p^a}5\r)F_{p^a}-\l(\f{p^a}5\r)(2F_{p^a}-F_{p^a}L_{p^a})
\\=&(L_{p^a}-1)\l(\l(\f{p^a}5\r)F_{p^a}-1\r)
\endalign$$
and hence it suffices to show that
$$\l(\f{p^a}5\r)F_{p^a}-1\eq\f{L_{p^a}-1}5\ (\mo\ p^2)\tag4.4$$
as $L_{p^a}\eq1\ (\mo\ p)$. Since $L_{p^a}\eq L_p\ (\mo\ p^2)$, and
$$F_{p^a}\eq\l(\f p5\r)^{a-1}F_p\ (\mo\ p^2)$$
by [S10, (2.5)], (4.4) is reduced to the case $a=1$. Note that
$$\l(\f p5\r)F_p-1=L_p-2F_{p-(\f p5)}-1\eq\f{L_p-1}5\ (\mo\ p^2)$$
since $L_p\eq1+\f52F_{p-(\f p5)}\ (\mo\ p^2)$ by the proof of [ST, Corollary 1.3].
The proof is now complete. \qed

\medskip
\noindent{\it Proof of Theorem 1.3}. Applying [ST, (2.4)] with $d=0$ and $n=p^a$, we get
$$\sum_{k=0}^{p^a-1}(-1)^k\bi{2k}k=\sum_{k=0}^{p^a-1}(-1)^k\bi{2p^a}kF_{2(p^a-k)}.\tag4.5$$
For each $k=1,\ldots,p^a-1$, clearly
$$\align(-1)^k\bi{2p^a}k=&(-1)^k\f{2p^a}k\bi{2p^a-1}{k-1}=-\f{2p^a}k\prod_{0<j<k}\l(1-\f{2p^a}j\r)
\\\eq&-\f{2p^a}k\(1-2p^a\sum_{0<j<k}\f1j\)\ (\mo\ p^3)
\endalign$$
and similarly
$$(-1)^k\bi{p^a}k\eq-\f{p^a}k\(1-p^a\sum_{0<j<k}\f1j\)\ (\mo\ p^3),$$
hence
$$(-1)^k\bi{2p^a}k\eq 4(-1)^k\bi{p^a}k+\f{2p^a}k\ (\mo\ p^3).$$
So, (4.5) yields that
$$\align&\sum_{k=0}^{p^a-1}(-1)^k\bi{2k}k-F_{2p^a}
\\\eq&4\sum_{k=1}^{p^a-1}\bi{p^a}k(-1)^kF_{2(p^a-k)}+2\sum_{k=1}^{p^a-1}\f{p^a}kF_{2(p^a-k)}
\\=&4\sum_{k=1}^{p^a-1}\bi{p^a}k(-1)^{p^a-k}F_{2k}+2\sum_{k=1}^{p^a-1}\f{p^a}{k}F_{2(p^a-k)}\ (\mo\ p^3).\endalign$$

Let $\al$ and $\beta$ be the two roots of the equation $x^2-x-1=0$. Then
$$\align\sum_{k=1}^{p^a-1}\bi{p^a}k(-1)^{p^a-k}F_{2k}
=&-F_{2p^a}+\sum_{k=0}^{p^a}\bi{p^a}k(-1)^{p^a-k}\f{\al^{2k}-\beta^{2k}}{\al-\beta}
\\=&\f{(\al^2-1)^{p^a}-(\beta^2-1)^{p^a}}{\al-\beta}-F_{2p^a}=F_{p^a}-F_{2p^a}.
\endalign$$
Thus, by the above and Proposition 4.2 we have
$$\align \sum_{k=0}^{p^a-1}(-1)^k\bi{2k}k\eq& F_{2p^a}+4(F_{p^a}-F_{2p^a})+2p^a\sum_{k=1}^{p^a-1}\f{F_{2(p^a-k)}}{k}
\\\eq&4F_{p^a}-3F_{2p^a}+2(F_{2p^a}-F_{p^a})+\l(\f{p^a}5\r)\f{(L_p-1)^2}5
\\=&2F_{p^a}-F_{2p^a}+\l(\f{p^a}5\r)\f{(L_p-1)^2}5\ (\mo\ p^3).
\endalign$$
Combining this with (4.3) we immediately get the desired (1.3).
This completes the proof of Theorem 1.3. \qed

\Ack. The work was supported by the National Natural Science
Foundation of China (Grant Nos. 10901078 and 11171140).
After we established Theorem 1.2,  we got Roberto Tauraso's generous permission (on Oct. 1, 2010)
to cite his conjecture (1.2) (cf. [T]). Thus
we would like to express our sincere thanks to Prof. R. Tauraso. We are also grateful to the two referees for their helpful comments.

\enddocument